\documentclass[11pt, reqno]{article}
\usepackage{amsmath, amssymb, amscd}
\usepackage{latexsym}

\newtheorem{thm}{Theorem}[section]
\newtheorem{cor}{Corollary}[section]
\newtheorem{pro}{Proposition}[section]

\newtheorem{lem}{Lemma}[section]

\newenvironment{prf}{\medskip 
\noindent {\em Proof.}}{\hfill $\square$ \\}
\newenvironment{rmk}{\medskip 
\noindent {\bf Remark.}}{\hfill \mbox{} \\}

\def\ms{\medskip}
\def\bs{\bigskip}
\def\nind{\noindent}

\def\ie{i.\ e.\ }

\def\ra{\rightarrow}
\def\({\left(}                   
\def\){\right)}                  
\def\[{\left[}                   
\def\]{\right]}                  
\def\lan{\langle}
\def\ran{\rangle}

\def\operatorname#1{{\rm #1\,}}                                
\def\op{\operatorname}

\def\tms{\times}

\def\grd{\nabla}         
\def\del{\partial}       
\def\lap{\triangle}

\def\a{\alpha}           
\def\b{\beta}          
\def\g{\gamma}           
\def\d{\delta}         
\def\e{\epsilon}         
\def\z{\zeta}            
\def\th{\theta}          
           
\def\lam{\lambda}        
\def\m{\mu}       
\def\n{\nu}       
\def\x{\xi}       
\def\p{\pi}       
\def\r{\rho}

\def\s{\sigma}    
\def\ph{\phi}     
     
\def\om{\omega}   
\def\G{\Gamma}

\def\Lam{\Lambda} 
       
\def\S{\Sigma}

\def\Om{\Omega}   

\def\R{\mathbb R}

\def\H{\mathbb H}

\def\bM{\overline M}
\def\bg{\overline g}
 
\numberwithin{equation}{section}
\title{On the $L^2$ Cohomology of a Convex Cocompact Hyperbolic Manifold}
\author{Xiaodong Wang \thanks{Department of Mathematics, MIT, Cambridge, 
MA 02139; \tt{Email:xwang@math.mit.edu}}}             
\date{September 2002}
\begin{document}
\maketitle
\begin{abstract}
We prove a vanishing theorem for a convex
cocompact hyperbolic manifold, which relates its $L^2$ cohomology 
and the Hausdorff
dimension of its limit set. The borderline case is shown to characterize the
manifold completely.
\end{abstract}
%%%%%%%%%%%%%%%%%%%%%%%%%%%%%%%%%%%%%%%%%%%%%%%%%%%%%%%%%%%%%%%%%%
%%%%%%%%%%%%%%%%%%%%%%%%%%%%%%%%%%%%%%%%%%%%%%%%%%%%%%%%%%%%%%%%%%
\section{Introduction}
The study of $L^2$ harmonic forms on a complete Riemannian manifold is
a very interesting and important subject.
In \cite{wa1} the author has studied $L^2$ harmonic 1-forms on a conformally
compact Einstein manifold and proved 
the following theorem.
\begin{thm}
Let $(M^{n+1}, g)$ be a conformally compact Einstein manifold.
\begin{enumerate}
\item If $\lam_0(g)>n-1$ then $\mathcal{H}^1(M)=0$. 
\item If $\lam_0(g)=n-1$ and $\mathcal{H}^1(M)\ne 0$, 
then $M$ is isometric to 
$\R\tms \S$, with warped product metric $dt^2+\cosh^2(t) h$, where $\S$
is compact and $h$ is a metric on $\S$ with $\op{Ric}(h)=-(n-1)h$.
\end{enumerate}
\end{thm}

Here $\lam_0(g)$ is the infimum of the $L^2$ spectrum of $-\lap$.
The proof hinges on the following inequality for a harmonic $1-$form
\begin{equation}\label{k1}
|\grd \th|^2\geq \frac{n+1}{n}|\grd |\th||^2
\end{equation}
and the characterization of the equality case. 

In this paper we use the same idea to study $L^2$ harmonic forms on a convex
cocompact hyperbolic manifold. There has been much work on this topic. 
We simply mention the paper by Mazzeo and Phillips \cite{mp} and the recent
work by Lott \cite{lo} and refer the reader to the reference therein for 
more background knowledge. Our main result is

\begin{thm}\label{main}
Let $M=\H^{n+1}/\G$ be a convex cocompact hyperbolic manifold
and $\d$ the Hausdorff dimension of the limit set of $\G$. Suppose
$\d>n/2$. Let $\mathcal{H}^p(M)$ be the space of 
$L^2$ harmonic $p-$forms.
\begin{enumerate}
\item If $p<n-\d$ then $\mathcal{H}^p(M)=0$. 
\item If $\d$ is an integer and $\mathcal{H}^{n-\d}(M)\ne 0$, 
then $M$ is a twisted warped product of $\H^{n-\d}$ and a compact hyperbolic
manifold of dimension $\d+1$ (described in detail in Section \ref{exp})
and $\op{dim}\mathcal{H}^{n-\d}(M)=1$.
\end{enumerate}
\end{thm}

Our proof is conceptually very simple. We use Bochner formula to prove
the vanishing theorem, but to get the sharp result we need a technical lemma
like (\ref{k1}).
It turns out that this 
inequality is an example of a refined Kato inequality and there are many
other examples in Riemannian geometry. Recently D. Calderbank, P. Gauduchon
and M. Herzlich \cite{cgh} have worked out a general principle which
covers all known examples and gives interesting 
new ones (T. Branson \cite{br} has
a different approach). As a special case their theorem implies that
\begin{equation}
|\grd \th|^2\geq \frac{n+2-p}{n+1-p}|\grd |\th||^2
\end{equation}
for a harmonic $p-$form on an $(n+1)$ dimensional Riemannian manifold.
Moreover the equality case is fully characterized. This result plays
a key role in the proof.
The proof
of the second part is a little bit involved and may have
some independent interest. In
this borderline case we have an $L^2$ harmonic form which satisfies
an overdetermined system of first order PDEs. The existence of such a
harmonic form gives rise to a splitting of $M$ and forces the metric
to be a twisted warped product.
It is surprising to have a situation where
the Bochner formula gives sharp results on
higher dimensional cohomology.

In closing the introduction we should mention 
the paper \cite{na} by Nayatani who proved a similar
result for compact Kleinian $n-$manifolds. His assumption requires
$\d<n/2-1$ while we assume that $\d> n/2$.
In some sense his result and ours are complementary. 

\bs

{\bf Acknowledgment:} I am indebted to Prof. Marc Herzlich
for drawing my attention to his joint paper with Calderbank and Gauduchon
\cite{cgh} which plays an important
role in the proof of the main theorem.
I wish to thank Professors Rick Schoen and 
Rafe Mazzeo for their constant encouragement and helpful discussions. 
Finally I want to thank the referee for many valuable comments and 
suggestions.

%%%%%%%%%%%%%%%%%%%%%%%%%%%%%%%%%%%%%%%%%%%%%%%%%%%%%%%%%%%%%%%%%%
\section{Preliminaries}\label{pre}
A complete hyperbolic manifold $(M^{n+1}, g)$ is the quotient of the 
unit ball $B^{n+1}$ by a torsion-free discrete group $\G$ of isometries of
$\H^{n+1}$. The limit
set $\Lam(\G)$ is defined to be the set of accumulation points in the 
sphere $S^n=\del B^{n+1}$ of an orbit
$\G(x)=\{\g(x)|\g\in\G\}$, where $x$ is a point in $B^{n+1}$.
$M$ is called geometrically finite if $\G$ has a fundamental domain bounded
by finitely many geodesic hyperplanes. $M$ is called convex cocompact
if the action of $\G$ on the hyperbolic convex hull of $\Lam(\G)$ in
$B^{n+1}$ has a compact fundamental region. Convex cocompact hyperbolic
manifolds can be characterized as geometrically finite hyperbolic 
manifolds without cusps. 

A convex cocompact hyperbolic manifold $M$ is conformally compact
in the sense that $\bM=M\sqcup (\Om(\G)/\G)$ is a compact manifold with
boundary and, if $r$ is a defining function (\ie a smooth function on
$\bM$ with first order zero on the boundary, positive on $M$), then
$\bg=r^2g$ extends as a regular metric on $\bM$. 
Its conformal infinity is the compact Kleinian manifold $\S=\Om(\G)/\G$.

The $L^2$ cohomology of a geometrically finite hyperbolic manifold was
studied by Mazzeo and Phillips \cite{mp}. We state their theorem for a convex
cocompact hyperbolic manifold.

\begin{thm}\label{mp}
Let $M=\H^{n+1}/\G$ be an orientable convex cocompact hyperbolic manifold.
There are natural isomorphisms
\begin{equation*}
\begin{cases}
\mathcal{H}^p\simeq H^p(M,\del M)& \text{if $p<(n+1)/2$} \\
\mathcal{H}^p\simeq H^p(M)& \text{if $p>(n+1)/2$}
\end{cases}
\end{equation*}
If $n+1$ is even then $\mathcal{H}^{(n+1)/2}$ is infinite dimensional.
\end{thm}

The asymptotics of such harmonic forms are also studied in detail in \cite{mp}.
To formulate the result we consider a neighborhood $\mathcal{U}$ of $\del M$
and use standard upper-half-space coordinates $(x,y)$ so that 
$\mathcal{U}\cap\del M=\{y=0\}$. Express $\om=\a+dy\wedge\b$, where $\a$
and $\b$ are a $p$ and $(p-1)$ form in $x$, respectively, depending 
parametrically on $y$.

\begin{thm}\label{mp2'}
Suppose $\om$ is an $L^2$ harmonic $p-$form in a neighborhood 
$\mathcal{U}$ of $\del M$.
Writing $\om=\a+dy\wedge\b$, the terms $\a,\b$ have complete asymptotic
expansions as $y\ra 0$, and in particular
\begin{align*}
\a&=\begin{cases}
   \a_{00}(x)y^{n-2p}+O(y^{n+1-2p}\log y),&  p<n/2 \\
   \a_{01}(x)y^2\log y+O(y^2),& p=n/2, 
   \end{cases} \\
\b&=\begin{cases}
   \b_{01}(x)y^{n+1-2p}\log y+O(y^{n+1-2p}),&  p<n/2 \\
   \b_{00}(x)y+O(y^2\log y),& p=n/2. 
   \end{cases}
\end{align*}
In all cases, the leading coefficients $\a_{00},\a_{01},\b_{00}, \b_{01}$
are $C^{\infty}$ and their estimates are uniform for $\om$'s bounded in 
$L^2(\mathcal{U})$. 
Similar expansions and rates of decay hold when $p\geq n/2+1$.
\end{thm}

The key to prove such asymptotic expansions is to construct a parametrix $Q$
for the Hodge Laplacian $\lap=dd^*+d^*d$ near the conformal infinity such that
$$Q\lap=I-R,$$
where $R$ is a smoothing operator. By construction the Schwartz kernel $R$ has
an asymptotic expansion, hence $\om$ has a similar expansion. In fact
it is shown in \cite{mp} that
\begin{equation}\label{sye}
\begin{split}
\a&\sim\sum^{\infty}_{j=0}\sum^{N_j}_{l=0}\a_{jl}(x)y^{n-2p+j}(\log y)^l, 
\quad N_0=0  \\
\b&\sim\sum^{\infty}_{j=0}\sum^{M_j}_{l=0}\b_{jl}(x)y^{n-2p+1+j}(\log y)^l,
\end{split}
\end{equation}
provided $p<n/2$. The forms $\a_{jl},\b_{jl}$ are $C^{\infty}$.
For $p=n/2$ there are similar results. For details
we refer to \cite{ma} and \cite{mp}. 

The above asymtotics given by Mazzeo and Phillips, 
while valid for solutions of $\lap \om = 0$
which are not necessarily closed and coclosed (e.g. as might
be needed when looking at solutions of $\lap \om= f$ where
f is compactly supported, hence zero in this boundary neighborhood),
can be improved if $\om$ is closed and coclosed. 
Write $\a=\sum_{|I|=p}\a_Idx^I$ and
$\b=\sum_{|J|=p-1}\b_Jdx^J$. From $d\om=0$ we get
\begin{equation*}
\sum_{|I|=p}\frac{\del \a_I}{\del y}dy\wedge dx^I
+\sum_{|I|=p,i}\frac{\del \a_I}{\del x^i}dx^i\wedge dx^I
-\sum_{|J|=p-1,j}\frac{\del \b_J}{\del x^j}dy\wedge dx^j\wedge dx^J=0.
\end{equation*}
Therefore
\begin{equation}
\frac{\del \a_I}{\del y}
=\sum_{J,j}\e^I_{jJ}\frac{\del \b_J}{\del x^j},
\end{equation}
where $\e^I_{jJ}=0$ unless $I=J\sqcup \{j\}$, in which case it is the
sign of the permutation $I \choose {jJ}$. This equation combined with
the asymptotic expansion (\ref{sye}) easily implies that the
coefficients $\a_{jl}=0$ in (\ref{sye}) for $j=0,1$. 
Similarly
we can prove that $\b_{0l}=0$ for $l\ne 0$ 
by using the equation $d^*\om=0$ combined
with the asymptotic expansion (\ref{sye}).

It is well-known that an $L^2$ harmonic form on a complete Riemannian 
manifold is both closed and coclosed. Therefore an $L^2$ harmonic form
on a convex cocompact hyperbolic manifold satisfies the  improved decay rate.
Though it is an elementary observation, 
we state it as a theorem for later reference. In this improved version
we do not need to formulate $p<n/2$ and $p=n/2$ separately.

\begin{thm}\label{mp2}
Suppose $\om$ is an $L^2$ harmonic $p-$form on a convex cocompact 
hyperbolic manifold $M$ with $p\leq n/2$. Then in a neighborhood 
$\mathcal{U}$ of $\del M$,
writing $\om=\a+dy\wedge\b$, the terms $\a,\b$ have complete asymptotic
expansions as $y\ra 0$, and in particular
\begin{align*}
\a&=O(y^{n+2-2p}\log y), \\ 
\b&=O(y^{n+1-2p}).
\end{align*}
Similar expansions and rates of decay hold when $p\geq n/2+1$.
\end{thm}

For a geometrically finite hyperbolic manifold, the asymptotics of an $L^2$
harmonic form at a cusp are also analyzed in \cite{mp}.

%%%%%%%%%%%%%%%%%%%%%%%%%%%%%%%%%%%%%%%%%%%%%%%%%%%%%%%%%%%%%%%%%%%%%%%%
\section{Special examples of convex cocompact \\ 
hyperbolic manifolds}\label{exp}
Let $(N, g_0)$ be a compact Riemannian manifold of dimension $k+1$ such that 
$\op{Ric}(g_0)=-kg_0$. Consider the following metric on $M=B^{n-k}\tms N$
$$g=\frac{4}{(1-|x|^2)^2}\left(dx^2+\frac{(1+|x|^2)^2}{4}g_0\right),$$
where $x$ is the coordinates on $B^{n-k}$. 
Then $g$ is a conformally compact 
Einstein metric. The conformal infinity is the $S^{n-k-1}\tms N$
with the product metric. If we use polar coordinates on the hyperbolic
space the metric can be written in the following form
\begin{equation}\label{xn}
g=dt^2+\sinh^2(t)d\z^2+\cosh^2(t)g_0, 
\end{equation}
where $d\z^2$ is the standard metric on ${S}^{n-k-1}$. It is obvious that
$(M,g)$ is the warped product of $\H^{n-k}$ and $(N,g_0)$.

If $(N, g_0)$ is hyperbolic, then $(M,g)$ is a convex cocompact
hyperbolic manifold. To see this we first 
consider the hyperbolic space $\H^{n+1}$
using the upper space model with coordinates $(r,x,y)$, where 
$r>0, x\in \R^k, y\in\R^{n-k}$. The hyperbolic metric is
\begin{equation*}
g=r^{-2}(dr^2+dx^2+dy^2).
\end{equation*} 
We introduce polar coordinates $y=\r\z$ on $\R^{n-k}$, 
with $\r>0,\z\in S^{n-k-1}$.
Then 
\begin{equation}
g=r^{-2}(dr^2+d\r^2+\r^2d\z^2+dx^2),
\end{equation} 
where $d\z^2$ is the standard metric on $S^{n-k-1}$.
We change coordinates by setting
\begin{equation}\label{cha}
r=s/\cosh (t), \quad \r=s \tanh (t).
\end{equation}
Straightforward calculation shows that in the new coordinates
\begin{equation}\label{warp}
g=dt^2+\sinh^2(t)d\z^2+\cosh^2(t) s^{-2}(ds^2+dx^2).
\end{equation}
This demonstrates that $\H^{n+1}$ is the warped product 
$\H^{n-k}\times\H^{k+1}$, since  
$$dt^2+\sinh^2(t)d\z^2$$ 
is exactly the hyperbolic metric in geodesic
polar coordinates,  and since 
$$s^{-2}(ds^2+dx^2)$$ 
is the hyperbolic metric on $\H^{k+1}$. This change of 
coordinates has a clear geometric meaning. $\H^{k+1}$ sits in $\H^{n+1}$
as the totally geodesic submanifold $\{y=0\}$. In our new coordinates
we simply view $\H^{n+1}$ by the exponential map on the normal bundle
of $\H^{k+1}$ in $\H^{n+1}$. It is easy to verify that $t$ as given
by (\ref{cha}) is the distance from the point $(r,x,y)$ to $\H^{k+1}$
with the closest point being $(s,x,0)$. Let $\G$ be the cocompact Kleinian
group such that $N=\H^{k+1}/\G$. There is a natural way to extend the action
of $\G$ to $\H^{n+1}$ ($n>k$) called the Poincar\'e extension. In terms of
the above description of $\H^{n+1}$ as the product $\H^{n-k}\times\H^{k+1}$
with the warped product metric (\ref{warp}), the extension is that $\G$
only acts on the second component. Hence $\H^{n+1}/\G=\H^{n-k}\tms N^{k+1}$ 
with the warped product metric given by (\ref{xn}).
\ms

There is a slightly more general construction. We give two equivalent
descriptions here.
Suppose that $E\ra N$ is a $\op{rank}=n-k$ flat $O(n-k)$ bundle. Such
a bundle is determined by its holonomy $\r:\G\ra O(n-k)$. We can cover
$N$ by open sets $\{U_{\a}\}$ with parallel trivialization 
$E|_{U_{\a}}\ra \R^{n-k}\tms U_{\a}$. Then on $E|_{U_{\a}}$ we can define
a hyperbolic metric using formula (\ref{xn}). Since the transition functions
$U_{\a}\cap U_{\b}\ra O(n-k)$ are locally constant, we get a global hyperbolic
metric which is apparently conformally compact. We call such a hyperbolic
manifold a twisted warped product. The second description we give is simpler.
We define a generalized Poincare extension of $\G-$action on $\H^{k+1}$
to $\H^{n+1}=\H^{n-k}\tms \H^{k+1}$ by 
using the homomorphism $\r:\G\ra O(n-k)$
$$\g\cdot (x,y)=(\r(\g)x,\g\cdot y).$$
Then $M=\H^{n+1}/\G$ is the twisted product. It is obvious that the limit
set of $M$ is a totally geodesic 
$S^k\subset S^n$. The convex core is the compact $k+1$ dimensional hyperbolic
manifold $N$ which is totally geodesic and the full manifold $M$ is
a vector bundle over $N$ with rank $n-k$.
Conversely we have the following
proposition whose proof is simple and hence omitted.
\begin{pro}
Let $M=\H^{n+1}/\G$ be a convex cocompact hyperbolic manifold. Suppose
the limit set is a totally geodesic
$S^k\subset S^n$, then $M$ is a twisted warped
product of $\H^{n-k}$ and a compact $k+1$ dimensional hyperbolic manifold.
\end{pro}

Suppose $k>n/2-1$ and both $N$ and $M$ are orientable,
we can describe $L^2$ harmonic $(n-k)-$forms on $M$ 
explicitly. First note $H^*(M)=H^*(N)$ for $N$ is a deformation retract
of $M$.
By Lefschetz duality and Mazzeo-Phillips theorem
\begin{equation*}
\mathcal{H}^{n-k}(M)\simeq H^{n-k}(M,\S)\simeq H^{k+1}(M)\simeq H^{k+1}(N)
\simeq \R.
\end{equation*}
If we introduce polar coordinates on the normal bundle of $N$ in
$M$, by the previous discussion the metric 
$g=dt^2+\sinh^2(t)d\z^2+\cosh^2(t)h$, where $h$ is the metric
on $N$ and $d\z^2$ is the standard metric on $S^{n-k-1}$. By calculation
one can show that all $L^2$ harmonic $(n-k)-$forms on $M$ are given by
the following formula
\begin{equation}
\om=c\frac{\sinh^{n-k-1}(t)}{\cosh^{k+1}(t)}dt\wedge \Theta,
\end{equation}
where $\Theta$ is the volume form on $S^{n-k-1}$ and $c$ is any constant.
It is easy to see that $|\om|=|c|\cosh^{-(k+1)}(t)$. 
Apparently the maximal level set is the convex core $N$. 
Note $t$ is the distance function to $N$.
%%%%%%%%%%%%%%%%%%%%%%%%%%%%%%%%%%%%%%%%%%%%%%%%%%%%%%%%%%%%%%%%%%%%
%%%%%%%%%%%%%%%%%%%%%%%%%%%%%%%%%%%%%%%%%%%%%%%%%%%%%%%%%%%%%%%%%%%%

\section{Proof of the main theorem}
We start to prove Theorem \ref{main}.
We prove part 1 by contradiction.
Suppose we have a nonzero $L^2$ harmonic form $\x$ of degree $p\leq n-\d$.
By Bochner formula 
\begin{equation}\label{bo}
0=(dd^*+d^* d)\x=\grd^*\grd\x+\mathcal{R}\x,
\end{equation}
where $\mathcal{R}\x=\th^k\wedge i_{e_l}R(e_k,e_l)\x$, if we choose
orthonormal frame $\{e_i\}$ for the tangent bundle and $\{\th^i\}$ the
dual frame for the cotangent bundle. The following lemma is well known.
For completeness we present the proof.

\begin{lem}\label{rterm}
Let $\om$ be a $p-$form on an $(n+1)-$dimensional manifold of constant
sectional curvature $\kappa$. Then
\begin{equation}
\mathcal{R}\om=p(n+1-p)\kappa\om.
\end{equation}
\end{lem}
\begin{prf}
As the metric has constant sectional curvature $\kappa$ we
have $$R(e_k,e_l)\th^i=-\kappa(\d_{li}\th^k-\d_{ki}\th^l).$$
Without loss of generality we assume 
$\om=\th^{1}\wedge\dotsb\wedge\th^{p}$.
We compute
\begin{align*}
\ &\th^k\wedge i_{e_l}R(e_k,e_l)\om \\
&=-\kappa\th^k\wedge i_{e_l}\(\sum^p_{i=1}(-1)^{i-1}
(\d_{li}\th^k-\d_{ki}\th^l)\wedge 
\th^{1}\wedge\dotsb\hat{\th^i}\dotsb\wedge\th^{p}\) \\
&=\kappa\sum^p_{i=1}(-1)^{i}\(\d_{li}\d^k_l-(n+1)\d_{ki}\)\th^k
  \wedge\th^{1}\wedge\dotsb\hat{\th^i}\dotsb\wedge\th^{p} \\
&\quad +\kappa\sum^p_{i=1}(-1)^{i}\d_{ki}\th^k\wedge\th^l 
  \wedge i_{e_l}\(\th^{1}\wedge\dotsb\hat{\th^i}\dotsb\wedge\th^{p}\) \\
&=np\kappa\om-p(p-1)\kappa\om \\
&=p(n+1-p)\kappa\om.
\end{align*}
\end{prf}

\nind
By this lemma and (\ref{bo}) we get
\begin{equation*}
\grd^*\grd\x=p(n+1-p)\x.
\end{equation*}
This easily implies 
\begin{equation}\label{lapn}
\frac{1}{2}\lap |\x|^2=|\grd\x|^2-p(n+1-p)|\x|^2.
\end{equation}
To proceed we need the following lemma which is 
Theorem 6.3.(ii) in \cite{cgh}, but
we state it in a way convenient for our purpose without introducing abstract
notations.

\begin{lem}\label{ka} 
Let $\x$ be a harmonic $p-$form (\ie $d\x=0$ and $\d\x=0$) on a Riemannian
manifold of dimension $n+1$, then
\begin{equation}
|\grd\x|^2\geq \frac{n+2-p}{n+1-p}|\grd |\x||^2.
\end{equation}
Moreover the equality holds iff there exists a $1-$form $\a$ with 
$\a\wedge\x=0$ such that
\begin{equation}
\grd \x=\a\otimes\x-\frac{1}{n+2-p}
        \sum^n_{i=0}\th^j\otimes(\th^j\wedge i_{\a^{\sharp}}\x),
\end{equation}
where $\{\th^0,\th^1,\ldots,\th^n\}$ is an orthonormal basis for the cotangent
bundle and $\a^{\sharp}$ is the vector dual to the $1-$form $\a$.
\end{lem} 
Let $f=|\x|$. By the above lemma, we get from (\ref{lapn}) 
$$\frac{1}{2}\lap f^2\geq \frac{n+2-p}{n+1-p}|\grd f|^2-p(n+1-p)f^2,$$
or, equivalently
$$f\lap f\geq \frac{1}{n+1-p}|\grd f|^2-p(n+1-p)f^2.$$
Let $\ph=f^{\b}$. We compute
\begin{align*}
\frac{1}{2}\lap \ph^2&=\ph\lap\ph+|\grd \ph|^2 \\
&=f^{\b}\[\b f^{\b-1}\lap f+\b(\b-1)f^{\b-2}|\grd f|^2\]+|\grd \ph|^2 \\
&=\b f^{2(\b-1)}\[f\lap f+(\b-1)|\grd f|^2\]+|\grd \ph|^2 \\
&\geq \b f^{2(\b-1)}\[\frac{1}{n+1-p}|\grd f|^2-p(n+1-p)f^2+(\b-1)|\grd f|^2\]
  +|\grd \ph|^2 \\
&=\(2\b-\frac{n-p}{n+1-p}\)\frac{1}{\b}|\grd\ph|^2-p(n+1-p)\b\ph^2.
\end{align*}
Let $\b=\frac{n-p}{n+1-p}$, then we have
\begin{equation}\label{fin}
\frac{1}{2}\lap \ph^2\geq |\grd\ph|^2-p(n-p)\ph^2.
\end{equation}
We now take a defining function $r$ such that near the conformal infinity
$\S$ the metric $g=r^{-2}(dr^2+h_r)$, where $h_r$ is an $r-$dependent family
of metrics on $\S$. Define $M^{\e}=\{x\in M|r(x)\geq\e\}$.
For $\e$ small enough this is a compact manifold with boundary. By (\ref{fin})
\begin{equation}\label{bdn}
\int_{M^{\e}}\(|\grd\ph|^2-p(n-p)\ph^2\)dV\leq 
\int_{\del M^{\e}}\ph\frac{\del\ph}{\del\n}d\s,
\end{equation}
where $\n$ is the outer unit normal of $\del M^{\e}$. By Theorem \ref{mp2}
we have
$$\ph=O(r^{n-p}).$$
Notice $\frac{\del\ph}{\del\n}$ is of the
same order as $\ph$.
Therefore we get
\begin{equation*}
\int_{\del M^{\e}}\ph\frac{\del\ph}{\del\n}d\s
=\int_{\S}O(\e^{n-2p})dV_{\e},
\end{equation*}
where $dV_{\e}$ is the volume form of $h_{\e}$ on $\S$.
Under the condition $n-p\geq \d> n/2$, 
the boundary term apparently goes to zero as 
$\e\ra 0$. Therefore
\begin{equation}
\int_{M}|\grd\ph|^2\leq p(n-p)\int_M \ph^2.
\end{equation}
According to Sullivan \cite{su} the infimum of the spectrum of $-\lap$
is given by $\lam_0=\d(n-\d)$.
If $n-p>\d$, then $\lam_0>p(n-p)$ and the above inequality is impossible. 
Hence $\mathcal{H}^p=0$. This finishes the proof of the first part.

Next, we prove the second part of Theorem \ref{main}.
Suppose that $\d$ is an integer and $\x$ is a nonzero $L^2$ harmonic
form of degree $p=n-\d$. Then from the previous discussion we must have
\begin{equation}\label{eig}
-\lap\ph=\d(n-\d)\ph,
\end{equation}
and, in view of Lemma \ref{ka}
\begin{equation}\label{grx}
\grd \x=\a\otimes\x-\frac{1}{\d+2}
        \sum^n_{i=0}\th^j\otimes(\th^j\wedge i_{\a^{\sharp}}\x)
\end{equation}
for a $1-$form $\a$ with $\a\wedge\x=0$. 

Before we plunge into the details, we describe our strategy. The existence
of a nontrivial solution of the overdetermined equations (\ref{grx}) will 
be used to show
that the regular level sets of $\ph$ are compact hypersurfaces of constant
mean curvature and they carry local splittings. 
Then by taking limit we prove that the maximum level set is a 
totally geodesic compact submanifold. The exponential map from its
normal bundle is then a diffeomorphism onto $M$ and gives the twisted
warped product structure. 
\ms

\nind
{\bf Step 1.}
By Harnack inequality the function
$\ph$ is everywhere positive on $M$. Let $c>0$ be a regular
value of $\ph$, then $\S_c=\ph^{-1}(c)$ is a compact hypersurface in $M$. 
Near a point $x\in \S_c$ we choose orthonormal basis of 1-forms 
$\{\th^0,\ldots,\th^n\}$ such that 
$\a=(\d+2)u\th^0$ where
$u$ is positive. As 
$\a\wedge\x=0$, we can write $\x=\th^0\wedge\om$ such that $\om$ contains
no components involving $\th^0$. By (\ref{grx}) we obtain the following
two equations
\begin{gather}
\grd_{e_0}(\th^0\wedge\om)=(\d+1)u\th^0\wedge\om, \label{g0} \\
\grd_{e_j}(\th^0\wedge\om)=-u\th^j\wedge\om, \quad j=1,\ldots,n.
\label{g2}
\end{gather}
Thus (recall $\ph=|\x|^{\d/(\d+1)}$)
\begin{gather}
e_0\ph=\frac{\d}{\d+1}|\x|^{-(\d+2)/(\d+1)}\lan\grd_{e_0}\x,\x\ran
=\d u\ph,\label{e0f} \\
e_j\ph=\frac{\d}{\d+1}|\x|^{-(\d+2)/(\d+1)}\lan\grd_{e_j}\x,\x\ran
=0, \quad j=1,\ldots,n.
\end{gather}
Therefore $e_0$ is the normal vector field of the hypersurface $\S_c$
and $e_1,\ldots, e_n$ are tangent to $\S_c$. Thus
\begin{equation}\label{e0}
\grd \ph=\d u\ph e_0.
\end{equation} 
We can write
\begin{equation}\label{sf}
\grd_{e_i}\th^0=\sum^n_{j=1}\Pi_{ij}\th^j, 
\end{equation}
where $\Pi_{ij}=\lan \grd_{e_i}e_0,e_j\ran$ is the second fundamental form
of the hypersurface $\S_c$.

Since $\grd_{e_0}(\th^0\wedge\om)=\grd_{e_0}\th^0\wedge\om
+\th^0\wedge\grd_{e_0}\om$ and $\grd_{e_0}\th^0$ has no $\th^0-$component,
we obtain from (\ref{g0})
\begin{gather}
\grd_{e_0}\th^0\wedge\om=0, \\
\th^0\wedge\(\grd_{e_0}\om-(\d+1)u\om\)=0.
\end{gather}
Similarly from (\ref{g2}) we get 
\begin{gather}
\th^0\wedge\grd_{e_j}\om=0, \label{par} \\
\(\grd_{e_j}\th^0+u\th^j\)\wedge\om=0.\label{htv}
\end{gather}
The equation (\ref{par}) implies that the tangential component of
$\grd_{e_j}\om$ is zero, \ie $\om$ restricted to $\S_c$ is parallel.

We introduce a distribution $E$ on $\S_c$ by defining
\begin{equation}
E_x=\{v\in T_x\S_c|v^*\wedge \om=0\}, \forall x\in \S_c.
\end{equation}
Let $E^{\perp}$ be the orthogonal complement of $E$. Then we have a 
decomposition 
\begin{equation}\label{dec}
T\S_c=E\oplus E^{\perp}.
\end{equation}
Both $E$ and $E^{\perp}$ are parallel for $\om$ is parallel. Obviously
$0\leq\op{rank}E\leq p-1=\op{deg}\om$.

The decomposition (\ref{dec}) gives a (local) splitting $\S_c=\S^1\times\S^2$
such that $g$ is the product of $g_1$ on $\S^1$ and $g_2$ on $\S_2$.
Hence $R_{\S_c}(X,Y,X,Y)=0, \forall X\in E, Y\in E^{\perp}$. By Gauss
equation
\begin{equation}\label{pia}
-|X|^2|Y|^2=-\Pi(X,X)\Pi(Y,Y)+\Pi(X,Y)^2.
\end{equation}
We can choose orthonormal bases on $E$ and $E^{\perp}$ such that in 
the corresponding coordinates
$$\Pi(X,X)=\sum_i\lam_i x_i^2, \quad \Pi(Y,Y)=\sum_j\m_j y_j^2.$$
By (\ref{pia}) we have
$$\Pi(X,Y)^2=-\sum_ix_i^2\sum_jy_j^2+\sum_i\lam_i x_i^2\sum_j\m_j y_j^2.$$
Fixing $x$, view both sides as quadratic forms in $y$. The right hand
side has no mixed terms $y_iy_j, i\ne j$. It follows that the linear
form $\Pi(X, Y)$ in $y$ involves only one of $y_j$'s. The same argument
works for $x$ while fixing $y$. Therefore by renumbering we can assume
$\Pi(X, Y)=ax_1y_1$. Then it is easy to see that $a=0$ and 
\begin{gather}
\Pi(X,Y)=0,\forall X\in E, Y\in E^{\perp} \\
\Pi(X,X)=|X|^2/\lam, \forall X\in E, \\
\Pi(Y,Y)=\lam|Y|^2, \forall Y\in E^{\perp},
\end{gather}
for some function $\lam$. We choose our basis such that 
$\th^1,\ldots,\th^s\in E^{\perp}$ and the rest in $E$, where 
$\d+1\leq s=\op{rank}E^{\perp}\leq n$.
By (\ref{sf}) (\ref{htv}) and the above three equations 
\begin{equation}\label{sft}
\Pi_{ij}=\begin{cases}
0,& i\ne j \\
-u, & i=j\leq s \\
-\frac{1}{u}, & i=j>s.
\end{cases}
\end{equation}
Hence the mean curvature of $\S_c$ is given by
\begin{equation}\label{mean}
H=-su-(n-s)/u.
\end{equation}
Again by Gauss equation and (\ref{sft})
\begin{align}
R_{\S^1}(X,Y,X,Y)=-1+1/u^2, 
\quad &\text{for orthonormal $X,Y\in T\S^1$} \\
R_{\S^2}(Z,W,Z,W)=-1+u^2, 
\quad &\text{for orthonormal $Z,W\in T\S^2$} 
\end{align}
It follows that $\op{Ric}_{\S_c}=s(-1+u^2)+(n-s)(-1+1/u^2)$. As
$\op{dim}\S_c=n\geq 3$, it is a standard consequence of the second 
Bianchi identity that $u$ is constant on $\S_c$.
In particular both $(\S^1,g_1)$ and $(\S^2, g_2)$ have constant sectional
curvatures. 
\ms

\nind
{\bf Step 2.}
Thus for $j=1,\ldots,n$
\begin{align*}
\lan \grd_{e_0}e_0,e_j\ran&=\lan \grd_{e_0}\(\grd \ph/|\grd \ph|\),e_j\ran \\
&=\lan \grd_{e_0}\grd \ph,e_j\ran/|\grd \ph| \\
&=\lan \grd_{e_j}\grd \ph,e_0\ran/|\grd \ph| \\
&=\frac{1}{|\grd \ph|}(e_je_0\ph-\grd_{e_j}e_0\ph) \\
&=\frac{1}{|\grd \ph|}\(\d e_j(u\ph)-\Pi_{ij}e_i\ph\) \\
&=0,
\end{align*}
where in the last step we use (\ref{e0f}) and the fact that $\ph$ and $u$  
are constant on $\S_c$. Therefore
\begin{equation}\label{geo}
\grd_{e_0}e_0=0.
\end{equation}

We claim that the constant $u$ is not $1$. Suppose that $u=1$, 
then $\S_c$ is 
a compact flat hypersurface in $M$. Its lifting in $\H^{n+1}$ is then
a horosphere which can be taken to be the hyperplane $S=\{y=a\}$ in the
upper space model for some $a>0$. 
Let $x$ and $\g\cdot x$ be in $S$,  which map to the same
point in the quotient $M=\H^{n+1}/\G$, where $\g\in\G$. Then $\g S=S$, and it
follows that $\g$ is a parabolic element. But this is impossible since $M$ has
no cusps.

Moreover $\op{rank}E=\op{deg}\om$, \ie $s=\d+1$. For otherwise we can write 
$\om=\th^{s+1}\dotsb\th^{n}\wedge\tau$ with $\tau$ a nontrivial parallel
form on $\S^2$. This would lead to a contradiction if we apply Lemma
\ref{rterm} on $\S^2$ which has nonzero constant curvature $-1+u^2$.

The equation (\ref{eig}) can be written as
\begin{equation}\label{laps}
-\d(n-\d)\ph=D^2\ph(e_0,e_0)+\lap_{\S_c}\ph+He_0\ph.
\end{equation}
The function $\ph$ being constant on $\S_c$, we get using (\ref{mean}) 
(\ref{dph})
\begin{align}\label{d2e0}
D^2\ph(e_0,e_0)&=-\d(n-\d)\ph-He_o\ph \\
&=-\d(n-\d)\ph+\((\d+1)+(n-\d-1)/u\)\d u\ph \\
&=\d(1+\d)u^2-\d\ph.
\end{align}
Combining these equations with (\ref{e0}) we obtain
\begin{equation}
D^2\ph(e_0,e_0)=\frac{1+\d}{\d\ph}|\grd\ph|^2-\d\ph.
\end{equation}
On the other hand $D^2\ph(e_0,e_i)=\lan \grd_{e_0}\grd\ph, e_i\ran
=|\grd\ph|\lan \grd_{e_0}e_0, e_i\ran=0$ for $i=1,\ldots, n$ while
\begin{align*}
D^2\ph(e_i,e_j)&=\lan\grd_{e_i}\grd\ph,e_j\ran \\
&=|\grd\ph|\lan \grd_{e_i}e_0,e_j\ran \\
&=\d u\ph\Pi_{ij}. \\
\end{align*}
Therefore we have
\begin{equation}
D^2\ph(e_i,e_j)=\begin{cases}
0,&  i\ne j \\
\frac{1+\d}{\d\ph}|\grd\ph|^2-\d\ph,& i=j=0 \\
-\frac{|\grd\ph|^2}{\d\ph},&  1\leq i=j\leq \d+1 \\
-\d\ph,& i=j>\d+1.
\end{cases}
\end{equation}
This shows that at any critical point of $\ph$ the Hessian $D^2\ph$ has 
constant rank $n+1-(\d+1)$. Therefore $N=\{x|\ph(x)=B\}$ is a nondegenerate
critical manifold of dimension $\d+1$,where $B=\op{max}\ph$. 
Let $h$ be the induced metric.
 
We show that $N$ is totally geodesic. Near $N$ we decompose $TM$ as the
direct sum of two subbundles according to the eigenspaces of $D^2\ph$.
We choose orthonormal basis $\{e_1,\ldots,e_{n+1}\}$ such that the first
$\d+1$ vector correspond to the eigenvalue $-\frac{|\grd\ph|^2}{\d\ph}$. 
On each regular level surface $\S_c$, 
$\{e_1,\ldots,e_{\d+1}\}$ span the distribution $E^{\perp}$ introduced
before. We know that $E$ and $E^{\perp}$ are parallel on $\S_c$, hence
$$\lan\grd_{e_i}e_k,e_j\ran=0$$
for $1\leq i,j\leq\d+1;k>\d+1$ and $\lan e_k,\grd\ph\ran=0$ while by
(\ref{sft})
$$\lan\grd_{e_i}\frac{\grd\ph}{|\grd\ph|},e_j\ran
=-u=-\frac{|\grd\ph|}{\d\ph}.$$ Therefore 
$$|\lan\grd_{e_i}X,e_j\ran|\leq \frac{|\grd\ph|}{\d\ph} $$
for any unit vector $X$ orthogonal to $e_1,\ldots,e_{\d+1}$. 
As $\grd\ph=0$ on $N$, we conclude that $N$ has zero second fundamental form
\ie totally geodesic. 
\ms

\nind
{\bf Step 3.}
We use the Ricci equation to compute the curvature
of the normal bundle $\mathcal{N}(N)$
\begin{align*}
\lan R^{\perp}_{VW}X,Y\ran&=R(V,W,X,Y)
+\sum^{\d+1}_{i=1}\(\Pi_X(V,e_i)\Pi_Y(W,e_i)-\Pi_X(W,e_i)\Pi_Y(V,e_i)\) \\
&=0.
\end{align*}
Hence the normal bundle is flat. Therefore we can choose our local orthonormal
frame on an open subset $U\subset N$ such that 
$e_{\d+2},\ldots, e_{n+1}$ are parallel sections of
$\mathcal{N}(N)$. 

Finally we consider the exponential map $\mathcal{N}(N)\ra M$ or locally
\begin{gather}
\psi:\R^+\tms S^{n-\d-1}\tms U\ra M,\\
\psi(t,\zeta,x)=exp_x\(t\sum^{n+1}_{i=d+2}\zeta^ie_i\).
\end{gather}
Given $V\in T_xN$ and $X\in T_{\zeta}S^{n-\d-1}$ we 
get Jacobi fields 
$V(t)=\psi_*(V)$ and $X(t)=\psi_*(X)$ along the geodesic 
$\g(t)=\psi(t,\zeta,x)$.
Note $V(0)=V$ and $\dot{V}(0)=\sum^{n+1}_{i=d+2}\zeta^i\grd_Xe_i=0$ 
because $e_{\d+2},\ldots, e_{n+1}$ are parallel sections of
$\mathcal{N}(N)$ and $N$ is totally geodesic. On the other hand
$X(0)=0, \dot{X}(0)=X$. Since the metric is hyperbolic, the Jacobi equation
is easy to solve to give
\begin{gather}
V(t)=\cosh (t)\mathcal{P}_t(V), 
X(t)=\sinh (t) \mathcal{P}_t(X),
\end{gather}
where $\mathcal{P}_t$ is the parallel translation from $T_xM$ to $T_{\g(t)}M$
along $\g$. Therefore in the geodesic polar coordinates $(t,\zeta,x)$ along
$N$ the metric takes the form
$$g=dt^2+\sinh^2(t) d\zeta^2+\cosh^2(t)h.$$
By (\ref{e0}) we have the following ODE on $\g$
\begin{equation}\label{dph}
\frac{d\ph}{dt}=-\d u\ph.
\end{equation}
We have negative sign here because $\dot{\g}=-e_0$ with our previous choice
of $e_0$.
We compute
\begin{align*}
\frac{d^2\ph}{dt^2}&=-\d\ph \frac{du}{dt}+\d u\frac{d\ph}{dt} \\
&=\d\ph \frac{du}{dt}+\d^2u^2\ph.
\end{align*}
On the other hand (\ref{d2e0}) gives us
\begin{equation*}
\frac{d^2\ph}{dt^2}=D^2\ph(e_0,e_0)=\d(1+\d)u^2-\d\ph.
\end{equation*}
Combining these two equations we obtain the ODE
\begin{equation}
\begin{cases}
\frac{du}{dt}=1-u^2& \\
u(0)=0.&
\end{cases}
\end{equation}
This can be easily solved and we get
\begin{gather}
u(t)=\frac{\sinh (t)}{\cosh (t)}, \\
\ph(t)=B\cosh^{-\d}(t).
\end{gather}
This shows that outside $\S$ the function $\ph$ is regular everywhere. 
Therefore $\psi:\mathcal{N}(N)\ra M$ is a diffeomorphism. This finishes
the proof.
\bs

\begin{rmk}
Part 1 of the theorem and its proof works for a geometrically 
finite hyperbolic manifold
whose only cusps are of maximum rank. A cusp of maximum rank is isometric
to $]1,\infty[\tms N$ with the metric $t^{-2}(dt^2+g_0)$, where $(N,g_0)$ is
a compact flat manifold. If we write $\om=\a+dt\wedge\b$, by \cite{mp} we
have
\begin{align*}
\a&=\begin{cases}
   \a_{0}(x)+O(e^{-\lam t}),&  k<n/2 \\
   O(e^{-\lam t}),& n/2\leq k\leq (n+1)/2, 
   \end{cases} \\
\b&=\begin{cases}
   \b_{0}(x)t+O(e^{-\lam t}),&  k<n/2 \\
   O(e^{-\lam t}),& n/2\leq k\leq (n+1)/2, 
   \end{cases}
\end{align*}
as $t\ra \infty$ for some $\lam>0$, where $\a_0$ and $\b_0$ are harmonic
forms on $N$.

When we do integration by parts on a compact domain in (\ref{bdn}), each
cusp gives rise to a boundary term 
$\int_{\{t\}\tms N}\ph\frac{\del\ph}{\del\n}d\s$ which tends to zero as 
$t\ra\infty$, by the asymptotics given above. Therefore the rest of the
argument goes without any change.

The asymptotics of $L^2$ harmonic forms near cusps of intermediate ranks
are also given in \cite{mp}, but the results are much more intricate. We do not
know whether the above proof can be generalized to cover the general case.
\end{rmk}

By Lefschetz duality and Theorem \ref{mp} we get the following corollary
from Theorem \ref{main}
\begin{cor}
Let $M=\H^{n+1}/\G$ be an orientable convex cocompact hyperbolic manifold
and $\d$ the Hausdorff dimension of the limit set of $\G$. Suppose
$\d> n/2$.
\begin{enumerate}
\item If $p>\d+1$ then $H^p(M,\R)=0$. 
\item If $\d$ is an integer and $H^{\d+1}(M,\R)\ne 0$, 
then $M$ is a twisted warped product of $\H^{n-\d}$ and a compact hyperbolic
manifold of dimension $\d+1$. 
\end{enumerate}
\end{cor}

\begin{rmk}
As shown by Izeki \cite{iz}, part one of the above corollary can be proven
by algebraic topology.
Let $\S=\Om(\G)/\G$ be the
conformal infinity which is a compact Kleinian $n-$manifold.
First by an idea in Schoen-Yau
\cite{sy} one can prove the relative homotopy groups $\p_i(M,\S)=0$ for 
$i<n-\d$. Then by Hurewicz isomorphsim theorem
$H^i(M,\S)=0$ for $i<n-\d$. By Lefschetz duality, $H^p(M,\R)=0$
for $p>\d+1$.

By Theorem \ref{mp}, this implies that
the $L^2$ cohomology is actually trivial if $\d<n/2$.
This can also be
easily seen from our approach, using (\ref{fin}) and the fact that 
$\lam_0=n^2/4$ when $\d\leq n/2$.   
By Mazzeo-Phillips theorem $\mathcal{H}^*(M)=0$ except for the middle 
dimension when $n+1$ is even. Therefore
the $L^2$ cohomology contains no useful information. However the interesting
work of Nayatani \cite{na}
shows that one can then
read off $\d$ from the cohomology of $\S$ when $\d<n/2-1$.
\end{rmk}

\end{document}